\documentclass[11pt,a4paper,reqno]{amsart}
\usepackage[T1]{fontenc}
\usepackage[utf8]{inputenc}
\usepackage{lmodern}
\usepackage[margin=28mm]{geometry}
\usepackage{amsmath,amssymb,mathtools}
\usepackage{microtype}
\usepackage{enumitem}
\usepackage[hidelinks]{hyperref}
\usepackage{bookmark}
\usepackage{doi}
\usepackage{xcolor}
\allowdisplaybreaks[2]
\numberwithin{equation}{section}
\newtheorem{theorem}{Theorem}[section]
\newtheorem{proposition}[theorem]{Proposition}
\newtheorem{lemma}[theorem]{Lemma}
\newtheorem{corollary}[theorem]{Corollary}
\theoremstyle{remark}
\newtheorem{remark}[theorem]{Remark}
\newcommand{\E}{\mathbb E}
\newcommand{\Prob}{\mathbb P}
\newcommand{\R}{\mathbb R}

\newcommand{\cK}{\mathcal K}

\newcommand{\cQ}{\mathcal Q}
\newcommand{\cA}{\mathcal A}
\newcommand{\cI}{\mathcal I}
\newcommand{\cF}{\mathcal F}
\newcommand{\bi}{\boldsymbol i}
\newcommand{\inj}{\mathrm{inj}}
\newcommand{\HS}{\mathrm{HS}}
\newcommand{\op}{\mathrm{op}}
\newcommand{\diam}{\operatorname{diam}}

\newcommand{\norm}[1]{\left\lVert #1\right\rVert}

\newcommand{\ip}[2]{\left\langle #1,#2\right\rangle}
\newcommand{\dd}{\,\mathrm d}
\newcommand{\ee}{\mathrm e}
\newcommand{\ext}{\mathfrak e}
\newcommand{\vf}{\mathfrak v}

\newcommand{\arxiv}[1]{\href{https://arxiv.org/abs/#1}{arXiv:#1}}

\hypersetup{pdftitle={Injective norms of inhomogeneous random tensors},
 pdfsubject={Inhomogeneous random tensors and injective norms},pdfauthor={Guozheng Dai and Ke Wang}}
\title[Inhomogeneous random tensors]
{Sharp Injective Norm Bounds for Inhomogeneous Gaussian Tensors and Sparse Extensions}
\author{Guozheng Dai}
\address{Department of Mathematics, Hong Kong University of Science and
Technology, Clear Water Bay, Kowloon, Hong Kong}
\email{guozhengdai@ust.hk}

\author{Ke Wang}
\address{Department of Mathematics, Hong Kong University of Science and
Technology, Clear Water Bay, Kowloon, Hong Kong}
\email{kewang@ust.hk}
\date{}
\makeatletter
\@ifundefined{subjclassname@2020}{%
  \@namedef{subjclassname@2020}{\textup{2020} Mathematics Subject Classification}%
}{}
\makeatother
\begin{document}
\begin{abstract}
\begingroup
We study the Euclidean injective norm of inhomogeneous random tensors
with independent entries. For Gaussian tensors, we establish optimal two-sided bounds for the expected norm, with constants
depending only on the tensor order. These bounds describe the interplay
between variance accumulation along fibres and the contribution of the
largest individual entries. In the matrix case, our result is consistent
with the characterizations in~\cite{LVHY18}.
We then extend the analysis to sparse Gaussian tensors with arbitrary
amplitudes and activation probabilities. The expected norm is characterized
by the expected largest absolute entry and the largest fibre standard
deviation after truncating amplitudes at this scale. Our proofs combine
geometric chaining with estimates for averaging over the random support.
\par\endgroup
\end{abstract}

\maketitle

\section{Introduction and main results}\label{sec:main}

Fix an integer $d\ge2$ and positive integers $n_1,\ldots,n_d$. Write
$[n]=\{1,\ldots,n\}$ and $\cI=[n_1]\times\cdots\times[n_d]$.
For a real tensor $U=(U_{\bi})_{\bi\in\cI}$, its \emph{Euclidean injective
norm} is defined by
\[
 \norm{U}_{\inj}
 =
 \sup_{\norm{x^{(r)}}_2\le1,\,1\le r\le d}
 \left|
   \sum_{\bi\in\cI}U_{\bi}\prod_{r=1}^d x^{(r)}_{i_r}
 \right|.
\]
Here $x^{(r)}\in\R^{n_r}$ for each $r$. Thus $\norm{U}_{\inj}$ is
the supremum of the absolute value of the associated multilinear
form over a product of Euclidean unit balls. When $d=2$, it
coincides with the usual matrix operator norm.

Estimating the injective norm of random tensors is important both
in theory and in applications. This problem extends spectral-norm
estimation for random matrices to higher-order random
arrays~\cite{Boed25} and also arises in the study of random
hypergraphs~\cite{ZZ21}, tensor principal component analysis with
random noise~\cite{MR14}, and tensor completion from randomly
sampled entries~\cite{JO14}, among other applications.
These connections motivate the search for sharp bounds that
capture how the distributions of the entries and their arrangement
affect the expected norm. In this paper, we seek sharp two-sided
bounds for the expected injective norm of random tensors with
independent entries.

Tomioka and Suzuki~\cite{TS14} obtained injective-norm upper bounds
of order $K\sqrt{(\sum_{r=1}^{d}n_r)\log(2d)}$ for
$n_1\times\cdots\times n_d$ tensors with independent centered
sub-Gaussian entries having a common sub-Gaussian parameter $K$.
For tensors with i.i.d.\ real or complex Gaussian entries,
Dartois and McKenna~\cite{DM24} used the Kac--Rice formula to obtain
high-probability upper bounds with refined leading constants.
For tensors with independent Bernoulli entries, Zhou and
Zhu~\cite{ZZ21,ZZ26} developed concentration and regularization
estimates allowing unequal success probabilities and deterministic
weights.
Lucca and Pesenti~\cite{LP26} established dimension-independent
bounds for independent entries with sufficiently small probabilities
of being nonzero, as well as Seginer-type two-sided estimates for
i.i.d.\ centered entries under an equivalence assumption on their
second and $(2+\varepsilon)$-th moments.
More recently, Lin and Lu~\cite{LL26} established a Seginer-type
fibre comparison for arbitrary centered integrable i.i.d.\ tensor
entries.

Another line of work studies random tensor series with independent
Gaussian or, more generally, sub-Gaussian
coefficients~\cite{BGJLR25,AA25}.
These results allow correlations between entries and also yield
upper bounds for the independent-entry model by taking the
deterministic summands to be coordinate tensors.
Specifically, for a deterministic real coefficient tensor
$A=(a_{\bi})_{\bi\in\cI}$, write
\[
 U=(a_{\bi}g_{\bi})_{\bi\in\cI}
 =\sum_{\bi\in\cI}g_{\bi}
   \bigl(a_{\bi}e_{i_1}\otimes\cdots\otimes e_{i_d}\bigr),
\]
where the $g_{\bi}$ are independent standard Gaussian variables.
Let $\cF_r$ denote the family of mode-$r$ fibres of $\cI$,
obtained by fixing all coordinates except $i_r$, and set
\[
 \sigma_r(A)=\max_{f\in\cF_r}
     \left(\sum_{\bi\in f}a_{\bi}^2\right)^{1/2},
 \qquad
 a_*=\max_{\bi\in\cI}|a_{\bi}|.
\]
After a symmetric embedding, the independent-entry estimate of
Bandeira et al.~\cite{BGJLR25} gives
\begin{equation}\label{eq:series-independent-comparison}
\begin{aligned}
 \E\norm{U}_{\inj}
 &\lesssim_d
 \left(\max_{1\le r\le d}\sigma_r(A)\right)
 \log\left(2+\sum_{r=1}^d n_r\right)\\
 &\quad+
 \max_{\substack{J\subseteq[d]\\ |J|\ge2}}
 a_*^{1-1/|J|}
 \max_{\bi_{J^c}}
 \left(\sum_{\bi_J}a_{\bi}^2\right)^{1/(2|J|)},
\end{aligned}
\end{equation}
where $\bi_J=(i_r)_{r\in J}$, the sums and maxima run over the
corresponding coordinate ranges, and the implicit constant depends
only on $d$.
Thus, this specialization involves both fibre variances and
variance sums over higher-dimensional slices.
Aden-Ali~\cite{AA25} removed the logarithmic multiplier on the
fibre term, while retaining the higher-dimensional slice terms.

Within the more restrictive independent-entry Gaussian model,
Boedihardjo~\cite{Boed25} obtained a bound that can be sharper than
the direct specializations of the general tensor-series estimates
of Bandeira et al.~\cite{BGJLR25} and Aden-Ali~\cite{AA25}.
His estimate involves only the fibre variances and the largest
coefficient.
When all mode dimensions are equal, $n_1=\cdots=n_d=n\ge2$,
\cite[Theorem~1.1]{Boed25} gives
\[
 \E\norm{U}_{\inj}
 \le
 \sqrt{2d}\sum_{r=1}^{d}\sigma_r(A)
 +Cd^3(\log n)^2a_*,
\]
where $C$ is a universal constant.
Here $\sigma_r(A)^2$ is the largest sum of entry variances along
a mode-$r$ fibre. Thus, the leading term reflects how variances
accumulate along fibres, while the correction term depends on
the largest coefficient.
This estimate also yields a bound of the form predicted by the
Euclidean type-2 conjecture discussed by
Lucca~\cite{Lucca25Conjecture}, with a polylogarithmic loss,
for Gaussian tensor series with independent entries.

The gain over the direct series estimates, and the remaining gap
to optimality, can both be seen in a diagonal example.
If $A$ has unit coefficients on the diagonal
$i_1=\cdots=i_d$ and is zero elsewhere, then
$\sigma_r(A)=a_*=1$ for every $r$.
For fixed $d$, the direct specializations of the bounds in
\cite{BGJLR25,AA25} contain a term of order
$n^{1/(2d)}$, whereas Boedihardjo's upper bound is of order
$(\log n)^2$.
Nevertheless, a logarithmic gap remains in this example, since
\[
 \E\norm{U}_{\inj}
 =
 \E\max_{1\le i\le n}|g_{i,\ldots,i}|
 \asymp \sqrt{\log n}.
\]
Thus, the correction $a_*(\log n)^2$ still overestimates the
actual scale of entrywise extremes.
This discrepancy motivates a finer characterization that captures
both variance accumulation along fibres and the actual scale of
entrywise extremes.
\subsection{Main results}
Our first main result gives a sharp two-sided characterization for
the independent-entry Gaussian model. It retains the fibre scale
and replaces the correction based on the largest coefficient by
an entrywise extreme term determined by the decreasing rearrangement
of the amplitudes. We then extend this characterization to sparse
Gaussian tensors with arbitrary amplitude and activation-probability
profiles.

We return to the rectangular index set
$\cI=[n_1]\times\cdots\times[n_d]$, with $d\ge2$ fixed.
Recall that $\sigma_r(A)$ denotes the largest Euclidean norm
of a mode-$r$ fibre of $A$. With
$\cF=\bigcup_{r=1}^d\cF_r$, set
\[
 \sigma(A)=\max_{1\le r\le d}\sigma_r(A)
 =\max_{f\in\cF}
   \left(\sum_{\bi\in f}a_{\bi}^2\right)^{1/2}.
\]
We write $x\wedge y=\min\{x,y\}$.
We write $X\lesssim_dY$ if $X\le C_dY$ for a positive constant
depending only on $d$, and $X\gtrsim_dY$ if $Y\lesssim_dX$.
We write $X\asymp_dY$ when both inequalities hold.
Omitted subscripts indicate absolute constants.

\begin{theorem}[Gaussian tensors with deterministic coefficients]\label{thm:gaussian}
Let $A=(a_{\bi})_{\bi\in\cI}$ be deterministic with $a_{\bi}\ge0$,
and let $(g_{\bi})_{\bi\in\cI}$ be independent standard Gaussian
variables.
Let $a_1^*\ge\cdots\ge a_{|\cI|}^*$ be the decreasing rearrangement of
the entries of $A$. Then
\begin{equation}\label{eq:unmasked}
 \E\norm{(a_{\bi}g_{\bi})_{\bi\in\cI}}_{\inj}
 \asymp_d
 \sigma(A)+\max_{1\le k\le|\cI|}a_k^*\sqrt{\log(k+1)}.
\end{equation}
\end{theorem}

The restriction to nonnegative amplitudes entails no loss for Gaussian
tensors: deterministic coefficient signs can be absorbed into the
independent Gaussian variables.

\begin{remark}[Comparison with earlier results]
\label{rem:comparison}
For each fixed order $d$, Theorem~\ref{thm:gaussian} gives a
sharp two-sided refinement of the Gaussian estimate of
Boedihardjo~\cite{Boed25}.
The fibre terms agree up to constants depending only on $d$, since
\[
 \sigma(A)\le\sum_{r=1}^d\sigma_r(A)\le d\,\sigma(A).
\]
In the case $n_1=\cdots=n_d=n\ge2$, the global correction
$a_1^*(\log n)^2$ is replaced by
\[
 \max_{1\le k\le n^d}a_k^*\sqrt{\log(k+1)}
 \le
 a_1^*\sqrt{d\log(n+1)}.
\]
Thus, for fixed $d$, the worst-case entrywise correction improves
from order $a_1^*(\log n)^2$ to order $a_1^*\sqrt{\log n}$,
while retaining the full decreasing amplitude profile can yield
further improvements.

The comparison constants in \eqref{eq:unmasked} depend on $d$,
and we do not optimize this dependence.
Boedihardjo's estimate, on the other hand, gives explicit dependence
on the tensor order through the factors $\sqrt{2d}$ and $d^3$.
Our refinement therefore concerns sharpness in the dimensions and
coefficient profiles for fixed $d$.

When $d=2$, Theorem~\ref{thm:gaussian} recovers, in an equivalent
deterministic form, the sharp operator-norm characterization of
Gaussian matrices due to Lata\l a, van Handel, and
Youssef~\cite[Corollary~4.1, with $p=\infty$]{LVHY18}.
\end{remark}

We next consider sparse Gaussian tensors, motivated by models in which
only a random subset of entries is observed or retained. Such random
support structures arise in covariance estimation with incomplete
observations, tensor completion and sparsification, randomized
embeddings, and random hypergraphs, among other applications.
For further background on these applications and related concentration
problems, see \cite{JO14,ZZ21,HWZ26,DHWZ25,DW26,DHWZ26,LP26}
and the references therein. Allowing both the activation probabilities
and the amplitudes to vary across the array accommodates nonuniform
sparsity and heterogeneous entry magnitudes.

Let $A=(a_{\bi})_{\bi\in\cI}$ and
$P=(p_{\bi})_{\bi\in\cI}$ be deterministic profiles with
$a_{\bi}\ge0$ and $p_{\bi}\in[0,1]$.
Let $\delta_{\bi}\sim\operatorname{Bernoulli}(p_{\bi})$ and
$g_{\bi}\sim N(0,1)$, with all these variables mutually independent,
and consider
\begin{equation}\label{eq:model}
 U=(a_{\bi}\delta_{\bi}g_{\bi})_{\bi\in\cI}.
\end{equation}
The amplitudes and activation probabilities may vary arbitrarily
across the array. 

The sparse characterization involves two deterministic quantities.
The first is the expected largest absolute entry of the sparse
Gaussian tensor in \eqref{eq:model}:
\begin{equation}\label{eq:expected-maximum}
 \ext(A,P)
 =\E_{\delta,g}\max_{\bi\in\cI}
       a_{\bi}\delta_{\bi}|g_{\bi}|
 =\E\max_{\bi\in\cI}|U_{\bi}|.
\end{equation}

The second quantity measures variance accumulation along fibres.
For $s\ge0$, define
\begin{equation}\label{eq:truncated-variance}
 v_s(A,P)
 =
 \max_{f\in\cF}
 \left(
   \sum_{\bi\in f}p_{\bi}(a_{\bi}\wedge s)^2
 \right)^{1/2},
 \qquad
 \vf(A,P)=v_{\ext(A,P)}(A,P),
\end{equation}
and put
\begingroup
\begin{equation}\label{eq:psi-sp}
 \Psi_{\mathrm{sp}}(A,P)
 =
 \ext(A,P)+\vf(A,P).
\end{equation}
\endgroup

\begin{theorem}[Sparse Gaussian tensors]\label{thm:sparse}
For every choice of dimensions and deterministic profiles $A,P$
as above, the tensor $U$ defined in \eqref{eq:model} satisfies
\begin{equation}\label{eq:main}
 \E\norm{U}_{\inj}
 \asymp_d\Psi_{\mathrm{sp}}(A,P).
\end{equation}
\end{theorem}

\begin{remark}[Sub-Gaussian amplitudes]\label{rem:subgaussian}
The upper bound in Theorem~\ref{thm:sparse} remains valid when
the Gaussian amplitudes are replaced by independent centered
sub-Gaussian variables.
More precisely, let $(\xi_{\bi})_{\bi\in\cI}$ be independent
of one another and of the Bernoulli masks, and suppose that
\[
 \E e^{t\xi_{\bi}}
 \le e^{K^2t^2/2},
 \qquad t\in\R,\quad \bi\in\cI,
\]
for some $K>0$.
Then the tensor
$U^{(\xi)}=(a_{\bi}\delta_{\bi}\xi_{\bi})_{\bi\in\cI}$
satisfies
\[
 \E\norm{U^{(\xi)}}_{\inj}
 \le C_dK\Psi_{\mathrm{sp}}(A,P).
\]

Indeed, conditional on the masks, the associated multilinear
process has sub-Gaussian increments controlled by $K$ times
the canonical metric of the corresponding Gaussian process.
The generic chaining upper bound and the Gaussian majorizing-measure
theorem~\cite{Talagrand14,Dirksen15} therefore give
\[
 \E_{\xi}
 \norm{(a_{\bi}\delta_{\bi}\xi_{\bi})_{\bi\in\cI}}_{\inj}
 \le
 CK\,\E_g
 \norm{(a_{\bi}\delta_{\bi}g_{\bi})_{\bi\in\cI}}_{\inj},
\]
where $C$ is an absolute constant.
Averaging over the masks and applying Theorem~\ref{thm:sparse}
yields the asserted upper bound.
\end{remark}

\begingroup
Bounding the extreme scale in terms of the largest amplitude and
expected support size gives the following simpler consequence of
Theorem~\ref{thm:sparse}.
\par\endgroup

\begingroup
\begin{corollary}[A variance-plus-maximum upper bound]
\label{cor:simple-upper}
For the sparse Gaussian tensor $U$ in \eqref{eq:model}, set
\[
 V(A,P)=\max_{f\in\cF}
 \Big(\sum_{\bi\in f}p_{\bi}a_{\bi}^2\Big)^{1/2},
 \qquad
 a_*=\max_{\bi\in\cI}a_{\bi},
 \qquad
 \Lambda=\sum_{\substack{\bi\in\cI\\a_{\bi}>0}}p_{\bi}.
\]
Then
\begin{equation}\label{eq:simple-upper}
 \E\norm{U}_{\inj}
 \lesssim_d V(A,P)+a_*\sqrt{\log(1+\Lambda)}.
\end{equation}
\end{corollary}

\par\endgroup

\begingroup
Here $V(A,P)^2$ is the largest total variance along a fibre,
and $\Lambda$ is the expected number of nonzero entries.
The bound avoids truncation and the ordered amplitude profile,
at the possible cost of sharpness. For an $n\times n$ Gaussian
matrix with $P\equiv1$, we have $\Lambda\le n^2$, so the bound becomes
$V(A,\mathbf1)+a_*\sqrt{\log(n+1)}$ up to universal constants.
This recovers the form of the Gaussian matrix bound of
Bandeira and van Handel~\cite{BvH16}.
\par\endgroup

\begingroup
Theorem~\ref{thm:sparse} also yields an upper bound for centered
Bernoulli tensors through symmetrization and Gaussian comparison.
\par\endgroup

\begin{corollary}[Bernoulli tensors]\label{cor:bernoulli-tensors}
Let $A=(a_{\bi})_{\bi\in\cI}$ be deterministic with $a_{\bi}\ge0$,
and let $\delta_{\bi}\sim\operatorname{Bernoulli}(p_{\bi})$ be
independent, where $p_{\bi}\in[0,1]$.
Define
\[
 U=(a_{\bi}\delta_{\bi})_{\bi\in\cI},
 \qquad
 q_{\bi}=2p_{\bi}(1-p_{\bi}),
 \qquad
 Q=(q_{\bi})_{\bi\in\cI}.
\]
Then
\begin{equation}\label{eq:bernoulli-fluctuation}
 \E\norm{U-\E U}_{\inj}
 \le C_d\Psi_{\mathrm{sp}}(A,Q),
\end{equation}
where $C_d$ depends only on $d$.
Here $\Psi_{\mathrm{sp}}(A,Q)$ is defined by
\eqref{eq:psi-sp} with $P$ replaced by $Q$.
\end{corollary}

\begingroup
\begin{remark}[Worst-case loss of the Bernoulli bound]
\label{rem:bernoulli-optimality}
The bound in Corollary~\ref{cor:bernoulli-tensors} can lose a factor
of order $\sqrt{\log n}$ on $[n]^d$, and this is the worst possible
loss for fixed $d$. More generally, with $N=|\cI|$,
\[
 \Psi_{\mathrm{sp}}(A,Q)
 \lesssim_d\sqrt{\log(N+1)}\,
 \E\norm{U-\E U}_{\inj}.
\]
This follows from Theorem~\ref{thm:sparse} and the Rademacher
contraction principle, which bounds the cost of replacing signs
by Gaussian variables by $\E\max_{\bi\in\cI}|g_{\bi}|
\lesssim\sqrt{\log(N+1)}$. Symmetrization changes the expected
norm by at most a factor of two. On $[n]^d$, we have $N=n^d$.

To see that this order is attained, let $A$ have unit coefficients
on a diagonal of length $n\ge2$ and vanish elsewhere, with activation
probability $1/2$ on the diagonal. Then
\[
 \norm{U-\E U}_{\inj}=\frac12
 \quad\text{almost surely},
 \qquad
 \Psi_{\mathrm{sp}}(A,Q)\asymp\sqrt{\log n}.
\]
Thus the factor $\sqrt{\log n}$ cannot be improved uniformly over
all profiles for fixed tensor order.
\end{remark}
\par\endgroup

\begingroup
Although the Bernoulli bound is not uniformly sharp, it can improve
existing estimates for inhomogeneous amplitudes. Zhou and
Zhu~\cite[Theorem~2.2]{ZZ26} obtain optimal expectation bounds at
every sparsity level for homogeneous tensors with $A\equiv1$ and
$p_{\bi}\equiv p\in(0,1/2]$ on $[n]^d$.
Their weighted extension also allows rectangular tensors and unequal
activation probabilities at every sparsity level~\cite[Corollary~2.3]{ZZ26}.
To compare the expectation bounds in a common regime, consider tensors
on $[n]^d$ with $p_{\bi}\equiv p\in(0,1/2]$ and $np\ge\log n$.
Their expectation upper bound then has scale $a_*\sqrt{np}$, whereas
\[
 \Psi_{\mathrm{sp}}(A,Q)
 \lesssim_d V(A,Q)+a_*\sqrt{\log n}
 \lesssim_d a_*\sqrt{np}.
\]
Our bound therefore improves this order whenever
$\log n=o(np)$ and $V(A,Q)=o(a_*\sqrt{np})$.
\par\endgroup

\subsection{Proof strategy}

Our argument separates the geometry of the tensor norm from the
averaging of the random support. We first prove
Theorem~\ref{thm:gaussian} for deterministic coefficients, and then
evaluate the random coefficient profiles obtained by conditioning
on the Bernoulli masks.

For the deterministic-coefficient estimate, we use the contraction
principle and uniform convexity estimates of van Handel~\cite{vH18},
and adapt the geometric-midpoint localization scheme proposed in
the candidate proof of Bednorz, Martynek, and
Meller~\cite[Section~6]{BMM26}.
The contribution here is to carry out this localization for a
multilinear parametrization of arbitrary fixed tensor order while
retaining the fibre variances and the ordered entrywise profile.
We use one fibre seminorm for each mode to control the variation
of the derivative. Their Gaussian averages are bounded by the
two terms in Theorem~\ref{thm:gaussian}, and the Hilbert--Schmidt
norm of the derivative is controlled by $\sigma(A)$.
These estimates allow square chaining in an auxiliary pseudometric
to yield the required chaining bound in the canonical pseudometric.
The transfer estimate is proved in full in
Subsection~\ref{sec:transfer}.

For the sparse model, the key contribution is to identify a common
scale that controls both the truncation error and the fluctuations
of the truncated fibre sums. The extreme parameter $\ext(A,P)$
provides this scale: amplitudes above a constant multiple of it
have a controlled total first-moment contribution, while truncation
gives a second-moment bound for the random entrywise profile.
Combined with symmetrization and an ordered maximum estimate,
these bounds compare $\E_\delta\sigma(A\odot\delta)$ with
$\vf(A,P)$ up to the additive term $\ext(A,P)$, as made precise
in Proposition~\ref{prop:averaging}. Here $\odot$ denotes entrywise
multiplication. The resulting comparison removes the remaining
expectation over the masks and explains the truncation appearing
in the deterministic characterization.

\section{Gaussian tensors with deterministic coefficients}
\label{sec:deterministic}

This section proves Theorem~\ref{thm:gaussian}. We retain the
fibre parameters $\sigma_r(A)$ and $\sigma(A)$ introduced in
Section~\ref{sec:main}. To simplify the notation in the proof,
for the deterministic nonnegative coefficient tensor
$A=(a_{\bi})_{\bi\in\cI}$ we set
\[
 \beta(A)=\max_{1\le j\le|\cI|}a_j^*\sqrt{\log(j+1)},
\]
where $(a_j^*)$ is the decreasing rearrangement of the
coefficients of $A$. Our goal is to prove
\begin{equation}\label{eq:coeff}
 \E\norm{(a_{\bi}g_{\bi})}_{\inj}
 \asymp_d\sigma(A)+\beta(A).
\end{equation}

The lower bound is straightforward and follows by testing
fibres and individual entries. The main difficulty lies in
the upper bound. We first establish the profile and chaining
estimates, then study the geometry of the multilinear
parametrization and prove the required localization bound.
The final subsection completes the proof of
Theorem~\ref{thm:gaussian}, including the lower bound.
Throughout the proof, we assume $A\ne0$, since the zero tensor
is immediate.

\subsection{Gaussian profiles and fibre seminorms}

The following classical estimate is recorded in
\cite[Lemma~4.5 and Remark~4.6]{vH17survey}. For completeness
and the reader's convenience, we include a short proof.

\begin{lemma}[Ordered Gaussian maximum]\label{lem:gmax}
Let $b_1,\ldots,b_M$ be nonnegative, and let
$b_1^*\ge\cdots\ge b_M^*$ denote their decreasing rearrangement.
For independent standard Gaussian variables $g_1,\ldots,g_M$,
\[
 \E\max_{j\le M}b_j|g_j|
 \asymp\max_{j\le M}b_j^*\sqrt{\log(j+1)}.
\]
\end{lemma}
\begin{proof}
The claim is immediate if all coefficients vanish.
By relabelling, we may assume $b_1\ge\cdots\ge b_M$.
The Gaussian tail bound and a union bound give, for $u\ge2$,
\[
 \Prob\left\{
 \max_{j\le M}b_j|g_j|
 >u\max_{k\le M}b_k\sqrt{\log(k+1)}
 \right\}
 \le2\sum_{j=1}^M(j+1)^{-u^2/2}
 \lesssim2^{-u^2/2}.
\]
Integration proves the upper bound. For the lower bound,
for every $1\le k\le M$,
\[
 \E\max_{j\le M}b_j|g_j|
 \ge b_k\,\E\max_{j\le k}|g_j|
 \gtrsim b_k\sqrt{\log(k+1)}.
\]
The last inequality follows from independence and the Gaussian
lower tail
\[
 \Prob\{|g|>t\}\gtrsim(1+t)^{-1}e^{-t^2/2},
 \qquad t\ge0.
\]
Taking the maximum over $k$ completes the proof.
\end{proof}

Recall that $\cF_r$ consists of the mode-$r$ fibres, obtained
by fixing all coordinates except the $r$th. Define
\begin{equation}\label{eq:fibrenorm}
 N_r(h)=\max_{f\in\cF_r}
 \left(\sum_{\bi\in f}a_{\bi}^2h_{i_r}^2\right)^{1/2},
 \qquad h\in\R^{n_r}.
\end{equation}
These seminorms are unconditional: if $|h_i|\le|k_i|$ for every
coordinate, then $N_r(h)\le N_r(k)$.

\begin{lemma}[Gaussian averages of the fibre seminorms]
\label{lem:auxGaussian}
For a standard Gaussian vector $g^{(r)}\in\R^{n_r}$,
\[
 \E N_r(g^{(r)})\lesssim\sigma_r(A)+\beta(A).
\]
\end{lemma}
\begin{proof}
Fibres on which all coefficients vanish may be omitted.
For each remaining $f\in\cF_r$, Jensen's inequality gives
\[
 \E\left(\sum_{\bi\in f}a_{\bi}^2
             (g^{(r)}_{i_r})^2\right)^{1/2}
 \le\left(\sum_{\bi\in f}a_{\bi}^2\right)^{1/2}
 \le\sigma_r(A).
\]
Moreover, the map
\[
 h\longmapsto
 \left(\sum_{\bi\in f}a_{\bi}^2h_{i_r}^2\right)^{1/2}
\]
is $\max_{\bi\in f}a_{\bi}$-Lipschitz in the Euclidean norm.
Since the fibres in $\cF_r$ are disjoint, their maximal
coefficients occupy distinct positions in $A$. Their decreasing
rearrangement is therefore bounded termwise by $(a_j^*)$.
Gaussian concentration \cite[Theorem~5.6]{BLM13} and a union
bound consequently give, for $u\ge2$,
\[
 \Prob\bigl\{N_r(g^{(r)})>\sigma_r(A)+u\beta(A)\bigr\}
 \le\sum_{j\ge1}(j+1)^{-u^2/2}
 \lesssim2^{-u^2/2}.
\]
Integration proves the claim. No independence between the fibre
norms is required.
\end{proof}

\subsection{Chaining and entropy estimates}\label{sec:prelim}

We follow the notation for entropy numbers and chaining functionals
in \cite[Section~2]{vH18}; see also \cite{Talagrand14} for background.
For a nonempty subset $U$ of a metric space $(X,d_0)$, define
\[
 \diam_{d_0}(U)=\sup_{x,y\in U}d_0(x,y)
\]
and, for each integer $j\ge0$,
\[
 e_j(U,d_0)
 =\inf_{\substack{F\subseteq X\\1\le|F|<2^{2^j}}}
   \sup_{x\in U}\inf_{y\in F}d_0(x,y).
\]
Thus $e_j(U,d_0)$ is the infimum of the radii of covers of $U$
by fewer than $2^{2^j}$ balls. The centres may lie in $X$;
requiring them to lie in $U$ changes this quantity by at most
a factor of two.

Let $S\subseteq X$ be nonempty. A finite partition $\cA_j$ of $S$
is a finite collection of pairwise disjoint nonempty subsets
whose union is $S$. Write $\cA_j(x)$ for its unique member
containing $x$. A sequence $(\cA_j)_{j\ge0}$ is
\emph{admissible} if
\[
 |\cA_j|<2^{2^j},\qquad
 \cA_{j+1}(x)\subseteq\cA_j(x)
 \qquad(j\ge0,\ x\in S).
\]
In particular, $\cA_0=\{S\}$. For $1\le p<\infty$, define
\[
 \gamma_{2,p}(S,d_0)
 =\inf_{\substack{(\cA_j)_{j\ge0}\\\text{admissible}}}
  \sup_{x\in S}
  \left(\sum_{j\ge0}
   2^{jp/2}\diam_{d_0}(\cA_j(x))^p\right)^{1/p},
\]
where the infimum runs over all admissible sequences of
partitions of $S$. We write $\gamma_2=\gamma_{2,1}$ for the
usual chaining functional; $\gamma_{2,2}$ uses a square sum
over the levels.
These definitions also apply to pseudometrics. When applying
metric-space results, we choose one representative from each
zero-distance class.

We use the notation
\[
 B_2^m=\{x\in\R^m:\norm{x}_2\le1\},\qquad
 B_{2,+}^m=B_2^m\cap[0,\infty)^m.
\]
For linear maps between Euclidean spaces, $\norm{\cdot}_{\HS}$
and $\norm{\cdot}_{\op}$ denote the Hilbert--Schmidt and
operator norms. When a seminorm $N$ appears in place of a
distance in $e_j$ or $\gamma_{2,p}$, the associated
pseudometric is $N(x-y)$.

We will use the following four estimates. The first controls
stochastic suprema, the second converts local entropy estimates
into a chaining bound, and the last two concern Euclidean balls
and their linear images.

\begin{proposition}[Standard chaining estimates]\label{prop:inputs}
The following estimates hold with absolute implicit constants.
\begin{enumerate}[label=(\roman*),leftmargin=*,itemsep=6pt]
\item\label{item:upperchaining}
Let $(Z_x)_{x\in S}$ be a separable centred process with
$Z_{x_0}=0$ almost surely for some $x_0\in S$. If
\[
 \E\exp\bigl(t(Z_x-Z_y)\bigr)
 \le\exp\left(\frac{t^2d_0(x,y)^2}{2}\right),
 \qquad t\in\R,\quad x,y\in S,
\]
then
\[
 \E\sup_{x\in S}|Z_x|\lesssim\gamma_2(S,d_0).
\]

\item\label{item:contraction}
Assume $\gamma_2(S,d_0)<\infty$, and let $\theta>0$ be a
sufficiently small absolute constant. If nonnegative functions
$s_j:S\to[0,\infty)$ satisfy
\[
 e_j(U,d_0)
 \le\theta\diam_{d_0}(U)+\sup_{x\in U}s_j(x)
\]
for every nonempty $U\subseteq S$ and every $j\ge0$, then
\[
 \gamma_2(S,d_0)
 \lesssim\sup_{x\in S}\sum_{j\ge0}2^{j/2}s_j(x).
\]

\item\label{item:convex}
For every seminorm $N$ on $\R^m$ and a standard Gaussian
vector $g\in\R^m$,
\[
 \gamma_{2,2}(B_2^m,N)\lesssim\E N(g).
\]

\item\label{item:ellipsoid}
For every linear map $L:\R^m\to\R^n$,
\[
 \sum_{j\ge0}2^j e_j(LB_2^m,\norm{\cdot}_2)^2
 \lesssim\norm L_{\HS}^2.
\]
\end{enumerate}
\end{proposition}

\begin{proof}
For part~\ref{item:upperchaining}, the increment assumption
allows us to apply the generic chaining upper bound
\cite{Talagrand14,Dirksen15}, which gives
\[
 \E\sup_{x\in S}|Z_x-Z_{x_0}|
 \lesssim\gamma_2(S,d_0).
\]
The conclusion follows from $Z_{x_0}=0$ almost surely.

For part~\ref{item:contraction}, the contraction principle
\cite[Theorem~3.1]{vH18}, applied with $\alpha=2$ and $p=1$,
yields
\[
 \gamma_2(S,d_0)
 \le C\theta\gamma_2(S,d_0)
 +C\sup_{x\in S}\sum_{j\ge0}2^{j/2}s_j(x)
\]
for an absolute constant $C$. Choosing $\theta\le(2C)^{-1}$
and using $\gamma_2(S,d_0)<\infty$, we may absorb the first
term on the right into the left-hand side to obtain
\[
 \gamma_2(S,d_0)
 \le2C\sup_{x\in S}\sum_{j\ge0}2^{j/2}s_j(x).
\]

For part~\ref{item:convex}, first suppose that $N$ is a norm.
The Euclidean unit ball is $2$-convex with an absolute
constant, as the parallelogram identity gives
\[
 \norm{\frac{x+y}{2}}_2
 \le1-\frac18\norm{x-y}_2^2,
 \qquad x,y\in B_2^m.
\]
The uniform convexity estimate
\cite[Theorem~5.8]{vH18}, with $\alpha=q=2$, and the dual
Sudakov inequality \cite[Lemma~8.3.6]{Talagrand14}
therefore give, respectively,
\[
 \gamma_{2,2}(B_2^m,N)
 \lesssim\sup_{j\ge0}2^{j/2}e_j(B_2^m,N)
 \lesssim\E N(g).
\]
For a seminorm, apply this estimate to the norm
$N+\varepsilon\norm{\cdot}_2$. Monotonicity gives
\[
 \gamma_{2,2}(B_2^m,N)
 \le\gamma_{2,2}(B_2^m,N+\varepsilon\norm{\cdot}_2)
 \lesssim\E N(g)+\varepsilon\E\norm{g}_2.
\]
The implicit constant is independent of $\varepsilon$,
so letting $\varepsilon\downarrow0$ proves the claim.

For part~\ref{item:ellipsoid}, the case $L=0$ is immediate.
By the singular value decomposition, $LB_2^m$ is a
possibly degenerate Euclidean ellipsoid whose semiaxes
are the nonzero singular values of $L$. The classical
entropy estimate for Hilbert-space ellipsoids
\cite[Section~2.5]{Talagrand14}, also used in
\cite[Section~7.3]{vH18}, bounds the weighted sum of squared
entropy numbers by the sum of squared semiaxes.
Since this latter sum equals $\norm{L}_{\HS}^2$, we obtain
\[
 \sum_{j\ge0}2^j e_j(LB_2^m,\norm{\cdot}_2)^2
 \lesssim\norm{L}_{\HS}^2.
\]
\end{proof}

We next combine the square-chaining estimates for the mode balls
on the parameter space used below. Identify
$\R^{n_1}\times\cdots\times\R^{n_d}$ with
$\R^{n_1+\cdots+n_d}$, equipped with
\[
 \norm{x}_2^2=\sum_{r=1}^d\norm{x^{(r)}}_2^2,
 \qquad x=(x^{(1)},\ldots,x^{(d)}).
\]
Set
\[
 \cK=B_{2,+}^{n_1}\times\cdots\times B_{2,+}^{n_d},
 \qquad
 \rho(x,y)=\sqrt d\sum_{r=1}^d
 N_r(x^{(r)}-y^{(r)}).
\]
Here $\rho$ is defined on the whole product space using the
fibre seminorms in \eqref{eq:fibrenorm}.

\begin{lemma}[Product square-chaining bound]\label{lem:product}
For the preceding parameter set and pseudometric,
\begin{equation}\label{eq:Pbound}
 \gamma_{2,2}(\cK,\rho)\lesssim_d\sigma(A)+\beta(A).
\end{equation}
\end{lemma}
\begin{proof}
For each $r$, choose admissible partitions
$(\cA_j^{(r)})_{j\ge0}$ of $B_2^{n_r}$, and set
$s=\lceil\log_2d\rceil+1$.
For $j\ge s$, take products of the cells of
$\cA_{j-s}^{(r)}$ restricted to $B_{2,+}^{n_r}$, omitting
empty cells; for $j<s$, use the trivial partition $\{\cK\}$.
The resulting partitions $(\cA_j)$ are nested and satisfy
\[
 |\cA_j|<2^{d2^{j-s}}\le2^{2^j},\qquad j\ge s,
\]
so they are admissible. For $x\in\cK$ and $j\ge s$, their
diameters satisfy
\[
 \diam_\rho(\cA_j(x))
 \le\sqrt d\sum_{r=1}^d
 \diam_{N_r}\bigl(\cA_{j-s}^{(r)}(x^{(r)})\bigr).
\]
Minkowski's inequality therefore gives
\[
 \begin{aligned}
 &\left(\sum_{j\ge s}2^j
          \diam_\rho(\cA_j(x))^2\right)^{1/2}\\
 &\qquad\le\sqrt d\,2^{s/2}\sum_{r=1}^d
 \left(\sum_{k\ge0}2^k
 \diam_{N_r}\bigl(\cA_k^{(r)}(x^{(r)})\bigr)^2
 \right)^{1/2}.
 \end{aligned}
\]
The initial levels are controlled by
\[
 \diam_\rho(\cK)
 \le\sqrt d\sum_{r=1}^d\diam_{N_r}(B_2^{n_r}),
\]
and these mode diameters already occur at level zero.
Taking the supremum over $x$ and then the infima over the
mode partitions yields
\[
 \gamma_{2,2}(\cK,\rho)
 \lesssim_d\sum_{r=1}^d\gamma_{2,2}(B_2^{n_r},N_r)
 \lesssim_d\sum_{r=1}^d\E N_r(g^{(r)}).
\]
The last inequality follows from
Proposition~\ref{prop:inputs}(\ref{item:convex}), and
Lemma~\ref{lem:auxGaussian} completes the proof.
\end{proof}

\subsection{Geometry of the multilinear parametrization}
\label{sec:geometry}

Define $\Phi:\R^{n_1+\cdots+n_d}\to\R^{\cI}$ by
\[
 \Phi(x)_{\bi}=a_{\bi}\prod_{r=1}^d x^{(r)}_{i_r},
\]
and equip $\cK$ with the canonical pseudometric
\[
 D(x,y)=\norm{\Phi(x)-\Phi(y)}_2
\]
of the Gaussian process $\ip{g}{\Phi(x)}$, where
$g=(g_{\bi})_{\bi\in\cI}$.
We write $\mathrm d\Phi(x)$ for the derivative of $\Phi$ at $x$.
The following estimates explain the role of the fibre seminorms:
the fibre variances control the size of the derivative, while
$\rho$ controls its variation and the error of linearization.

\begin{lemma}[Derivative estimates]\label{lem:tangent}
For $x,y\in\cK$,
\begin{align}
 \norm{\mathrm d\Phi(x)}_{\HS}
 &\le\sqrt d\,\sigma(A),\label{eq:Rbound}\\
 \norm{\mathrm d\Phi(x)-\mathrm d\Phi(y)}_{\op}
 &\le\rho(x,y).\notag
\end{align}
Moreover, for $x,y,z\in\cK$,
\begin{equation}\label{eq:linearization}
 \norm{\Phi(x)-\Phi(y)-\mathrm d\Phi(z)(x-y)}_2
 \le[\rho(x,z)+\rho(y,z)]\norm{x-y}_2.
\end{equation}
\end{lemma}
\begin{proof}
For $h=(h^{(1)},\ldots,h^{(d)})$, differentiation gives
\[
 \bigl(\mathrm d\Phi(x)h\bigr)_{\bi}
 =a_{\bi}\sum_{r=1}^d h^{(r)}_{i_r}
       \prod_{s\ne r}x^{(s)}_{i_s}.
\]
Summing the squared lengths of the columns of the derivative,
and then summing over the mode-$r$ fibres, gives
\[
 \begin{aligned}
 \norm{\mathrm d\Phi(x)}_{\HS}^2
 &=\sum_{r=1}^d\sum_{\bi\in\cI}a_{\bi}^2
       \prod_{s\ne r}(x^{(s)}_{i_s})^2\\
 &\le\sum_{r=1}^d\sigma_r(A)^2
       \prod_{s\ne r}\norm{x^{(s)}}_2^2
 \le d\sigma(A)^2.
 \end{aligned}
\]

For the derivative difference, fix $r$ and telescope
$\prod_{t\ne r}x^{(t)}_{i_t}-\prod_{t\ne r}y^{(t)}_{i_t}$
one mode at a time. In the term where mode $s\ne r$ changes,
the remaining factors come from the corresponding modes of
$x$ or $y$. Summing the squared coordinates first over the
mode-$s$ fibres, using \eqref{eq:fibrenorm}, bounds the norm
of this term by
\[
 N_s(x^{(s)}-y^{(s)})\norm{h^{(r)}}_2,
\]
since every remaining mode vector has Euclidean norm at most one.
The triangle inequality and Cauchy--Schwarz now yield
\[
 \begin{aligned}
 \norm{\bigl(\mathrm d\Phi(x)-\mathrm d\Phi(y)\bigr)h}_2
 &\le\sum_{r=1}^d\norm{h^{(r)}}_2
       \sum_{s\ne r}N_s(x^{(s)}-y^{(s)})\\
 &\le\sqrt d\,\norm h_2
       \sum_{s=1}^dN_s(x^{(s)}-y^{(s)})
 =\rho(x,y)\norm h_2.
 \end{aligned}
\]
Finally, $\cK$ is convex and
\[
 \Phi(x)-\Phi(y)-\mathrm d\Phi(z)(x-y)
 =\int_0^1
 \bigl[\mathrm d\Phi((1-t)y+tx)-\mathrm d\Phi(z)\bigr]
 (x-y)\,\dd t.
\]
Since $\rho$ is induced by a seminorm,
\[
 \rho((1-t)y+tx,z)
 \le(1-t)\rho(y,z)+t\rho(x,z).
\]
Integration proves \eqref{eq:linearization}.
\end{proof}

The localization argument also uses the coordinatewise geometric
midpoint from \cite[Section~6.1]{BMM26}. Its compatibility with
$\Phi$ gives the following squared-norm identity and distance bounds.

\begin{lemma}[Geometric midpoint]\label{lem:midpoint}
For $x,y\in\cK$, define
\[
 m(x,y)^{(r)}_i=\sqrt{x^{(r)}_iy^{(r)}_i}.
\]
Then $m(x,y)\in\cK$ and
\begin{align}
 \norm{m(x,y)}_2^2
 &=\frac{\norm{x}_2^2+\norm{y}_2^2-\norm{x-y}_2^2}{2},\notag\\
 D(x,m(x,y))&\le D(x,y),\label{eq:midD}\\
 \rho(x,m(x,y))&\le\rho(x,y).\label{eq:midrho}
\end{align}
\end{lemma}
\begin{proof}
For each $r$, Cauchy--Schwarz gives
\[
 \norm{m(x,y)^{(r)}}_2^2
 =\ip{x^{(r)}}{y^{(r)}}
 \le\norm{x^{(r)}}_2\norm{y^{(r)}}_2\le1.
\]
Thus $m(x,y)\in\cK$. Summing over $r$ and applying polarization
proves the squared-norm identity.
Since the coefficients and coordinates are nonnegative,
\[
 \Phi(m(x,y))_{\bi}
 =\sqrt{\Phi(x)_{\bi}\Phi(y)_{\bi}}.
\]
The scalar inequality
\[
 |u-\sqrt{uv}|\le|u-v|,\qquad u,v\ge0,
\]
applied to the coordinates of $\Phi(x)$ and $\Phi(y)$ gives
\eqref{eq:midD}. Applying it to the mode coordinates and using
the unconditionality of $N_r$ gives \eqref{eq:midrho}.
\end{proof}

\subsection{Localization and the canonical pseudometric}
\label{sec:transfer}

We now combine the square-chaining estimate, the derivative bounds,
and the geometric midpoint identity. The argument adapts the
localization method of \cite[Section~6]{BMM26} to $\Phi$.

\begin{proposition}[Chaining bound for the canonical pseudometric]
\label{prop:transfer}
With the preceding definitions,
\[
 \gamma_2(\cK,D)\lesssim_d\sigma(A)+\beta(A).
\]
\end{proposition}
\begin{proof}
We apply the contraction principle in
Proposition~\ref{prop:inputs}(\ref{item:contraction}).
To verify its entropy condition, we will cover each subset of
$\cK$ by first approximating its points through a regularized
minimization problem and then covering the resulting set by
linearization. The approximation errors and the costs of the
linear covers will give the functions $s_j$.
Fix a sufficiently small absolute $\theta>0$, and let $b>0$
be chosen at the end of the proof.

For $t>0$ and $z\in\cK$, set
\[
 J_t(z)=\min_{w\in\cK}
 \bigl\{\norm w_2^2+tD(z,w)+t^2\rho(z,w)^2\bigr\},
\]
and choose a minimizer $\pi_tz$, which exists by compactness
and continuity. The two distance penalties control the errors
in the canonical pseudometric $D$ and in the auxiliary
pseudometric $\rho$ used for linearization.
Accordingly, put $t_j=b2^{j/2}$ for $j\ge-1$
and, for $j\ge0$, write
\[
 \epsilon_j(z)=D(z,\pi_{t_j}z),\qquad
 \eta_j(z)=\rho(z,\pi_{t_j}z).
\]
Testing $w=z$ gives $0\le J_t(z)\le\norm z_2^2\le d$.
Testing $\pi_{t_j}z$ in the minimization problem at scale
$t_{j-1}$ gives
\[
 J_{t_j}(z)-J_{t_{j-1}}(z)
 \ge(1-2^{-1/2})t_j\epsilon_j(z)
     +\frac12t_j^2\eta_j(z)^2.
\]
The left-hand side telescopes, so
\begin{equation}\label{eq:budget}
 \sup_{z\in\cK}\sum_{j\ge0}
 \bigl[t_j\epsilon_j(z)+t_j^2\eta_j(z)^2\bigr]
 \lesssim d.
\end{equation}

To turn these approximations into covers, we need to control
both the variation of $\mathrm d\Phi$ and the size of its linear
images. For the first purpose, the linearization estimate
\eqref{eq:linearization} suggests working inside sets of
controlled $\rho$-diameter. By \eqref{eq:Pbound}, choose
admissible partitions $(\cQ_j)_{j\ge0}$ of $\cK$ such that
\begin{equation}\label{eq:partitionbudget}
 \sup_{z\in\cK}\sum_{j\ge0}2^j
 \diam_\rho(\cQ_j(z))^2
 \lesssim_d[\sigma(A)+\beta(A)]^2.
\end{equation}
For the second purpose, record the covering radii of the
derivative images by setting
\[
 \lambda_j(z)=e_j\bigl(
 \mathrm d\Phi(z)B_2^{n_1+\cdots+n_d},\norm{\cdot}_2\bigr).
\]
Proposition~\ref{prop:inputs}(\ref{item:ellipsoid}) and
\eqref{eq:Rbound} yield
\begin{equation}\label{eq:lambdabudget}
 \sup_{z\in\cK}\sum_{j\ge0}2^j\lambda_j(z)^2
 \lesssim_d\sigma(A)^2.
\end{equation}
Thus the partition diameters and the derivative covering radii
both have the square-summability needed below.

We now construct the cover required by the contraction
principle. Fix $j\ge1$ and a nonempty $U\subseteq\cK$.
We use the partition $\cQ_{j-1}$ and a level-$(j-1)$ cover
inside each of its cells. This shift allocates the available
points between the two operations: fewer than $2^{2^{j-1}}$
cells, each covered by fewer than $2^{2^{j-1}}$ points, give
fewer than $2^{2^j}$ points in total, as required for $e_j(U,D)$.
Fix a cell $Q\in\cQ_{j-1}$ with $U_0=U\cap Q\ne\varnothing$,
and denote the Euclidean diameter of its approximating points by
\[
 h=\diam_2\bigl(\pi_{t_j}(U_0)\bigr).
\]

The minimizing property of $\pi_{t_j}$ makes this diameter
small enough for linearization to be useful. Indeed, take
$z,z'\in U_0$, set $x=\pi_{t_j}z$ and $y=\pi_{t_j}z'$, and
exchange the two pairs if necessary so that
$\norm x_2^2\ge\norm y_2^2$.
Using the geometric midpoint from Lemma~\ref{lem:midpoint}
as a competitor in the minimization defining $x$, we obtain
\[
 \begin{aligned}
 \frac12\norm{x-y}_2^2
 &\le\norm x_2^2-\norm{m(x,y)}_2^2\\
 &\le t_j\bigl[D(z,m(x,y))-D(z,x)\bigr]\\
 &\qquad+t_j^2\bigl[\rho(z,m(x,y))^2-\rho(z,x)^2\bigr]\\
 &\le t_jD(x,y)+t_j^2\rho(z,m(x,y))^2.
 \end{aligned}
\]
Here the last inequality uses \eqref{eq:midD}.
Since $z,z'\in Q$, the triangle inequality and
\eqref{eq:midrho} give
\[
 \begin{aligned}
 D(x,y)&\le\diam_D(U)+2\sup_{z\in U_0}\epsilon_j(z),\\
 \rho(z,m(x,y))
 &\le\diam_\rho(Q)+3\sup_{z\in U_0}\eta_j(z).
 \end{aligned}
\]
Taking the supremum over $z,z'\in U_0$ therefore yields
\begin{equation}\label{eq:localdiam}
 \begin{aligned}
 h^2\lesssim{}&
 t_j\left[\diam_D(U)+\sup_{z\in U_0}\epsilon_j(z)\right]\\
 &+t_j^2\left[\diam_\rho(Q)^2+
                   \sup_{z\in U_0}\eta_j(z)^2\right].
 \end{aligned}
\end{equation}

Fix any $z_0\in U_0$. The inclusion
\[
 \mathrm d\Phi(z_0)
 \bigl(\pi_{t_j}(U_0)-\pi_{t_j}z_0\bigr)
 \subseteq h\,\mathrm d\Phi(z_0)B_2^{n_1+\cdots+n_d}
\]
allows us to cover the linear image by fewer than
$2^{2^{j-1}}$ balls with radii arbitrarily close to
$h\lambda_{j-1}(z_0)$. Assign each point of $U_0$ to a ball
containing the linear image of its approximation, and choose
one original point from each nonempty group as a covering centre.
If $x,y\in\pi_{t_j}(U_0)$ correspond to the same group, then,
up to an arbitrarily small error,
\[
 \begin{aligned}
 \norm{\mathrm d\Phi(z_0)(x-y)}_2
 &\le2h\lambda_{j-1}(z_0),\\
 \rho(x,z_0),\ \rho(y,z_0)
 &\le\diam_\rho(Q)+\sup_{z\in U_0}\eta_j(z).
 \end{aligned}
\]
Apply \eqref{eq:linearization}, use $\norm{x-y}_2\le h$,
and add the two $D$-approximation errors to return to the
original points. Letting the error in the covering radii
tend to zero gives
\begin{equation}\label{eq:rawcover}
 \begin{aligned}
 e_{j-1}(U_0,D)\lesssim{}&\sup_{z\in U_0}\epsilon_j(z)\\
 &+h\left[\lambda_{j-1}(z_0)+\diam_\rho(Q)
                      +\sup_{z\in U_0}\eta_j(z)\right].
 \end{aligned}
\end{equation}
The centres lie in $U_0$. If $h=0$, all approximating points
coincide, and one point of $U_0$ covers it with radius at most
$2\sup_{z\in U_0}\epsilon_j(z)$, so the same estimate holds.

Substitute \eqref{eq:localdiam} into \eqref{eq:rawcover} and
apply Young's inequality, with coefficient $\theta$ for the
term involving $\diam_D(U)$.
For the fixed $\theta$, the resulting covering radius is at most
\[
 \begin{aligned}
 \theta\diam_D(U)+C_\theta\Bigl\{
 &\sup_{z\in U_0}\epsilon_j(z)\\
 &+t_j\bigl[\lambda_{j-1}(z_0)^2+\diam_\rho(Q)^2
                    +\sup_{z\in U_0}\eta_j(z)^2\bigr]\Bigr\}.
 \end{aligned}
\]
This identifies the functions needed in the contraction
principle: for $j\ge1$, define
\begin{equation}\label{eq:sj}
 s_j(z)=C_\theta\left[
 \epsilon_j(z)+t_j\left(
 \lambda_{j-1}(z)^2+
 \diam_\rho(\cQ_{j-1}(z))^2+\eta_j(z)^2
 \right)\right],
\end{equation}
with $C_\theta$ enlarged if necessary.
On $U_0$, the partition term is constant because
$\cQ_{j-1}(z)=Q$, and the derivative term at $z_0$ is
bounded by its supremum over $U_0$.
Each of the four nonnegative terms in the preceding radius
bound is thus controlled by the supremum of the bracket in
\eqref{eq:sj}; their sum costs only an additional absolute
factor. Consequently, $U_0$ has a cover of the stated cardinality
with radius arbitrarily close to
\[
 \theta\diam_D(U)+\sup_{z\in U_0}s_j(z).
\]
Combining these covers over the cells of $\cQ_{j-1}$ gives
fewer than $2^{2^j}$ centres in $U$, and hence
\[
 e_j(U,D)\le\theta\diam_D(U)+\sup_{z\in U}s_j(z).
\]
For $j=0$, the same inequality holds with
$s_0(z)=\diam_D(\cK)$.

It remains to sum the covering costs and apply the contraction
principle. Since $\sup_{z\in\cK}\norm{\Phi(z)}_2
\le\max_{\bi}a_{\bi}$, the initial level contributes at most
$\diam_D(\cK)\le2\sigma(A)$.
For the remaining levels, the identities
\[
 2^{j/2}\epsilon_j(z)=\frac{t_j\epsilon_j(z)}b,
 \qquad
 2^{j/2}t_j\eta_j(z)^2=\frac{t_j^2\eta_j(z)^2}b
\]
show that \eqref{eq:budget} controls the two approximation
errors by a constant multiple of $d/b$.
The other two terms in \eqref{eq:sj} are controlled by
\eqref{eq:partitionbudget} and \eqref{eq:lambdabudget}, because
$2^{j/2}t_j=b2^j$ and the shift from $j$ to $j-1$ changes
each square sum only by a factor of two. Therefore
\[
 \sup_{z\in\cK}\sum_{j\ge0}2^{j/2}s_j(z)
 \lesssim_d\sigma(A)+\frac1b+b[\sigma(A)+\beta(A)]^2.
\]

For completeness, \eqref{eq:Rbound} and convexity of $\cK$
give $D(x,y)\le\sqrt d\,\sigma(A)\norm{x-y}_2$.
Finite-dimensional Euclidean covering estimates then imply
$\gamma_2(\cK,D)<\infty$, as required by the contraction
principle. If $D$ is only a pseudometric, apply that principle
to a set containing one representative of each zero-distance
class, restrict the $\rho$-partitions to this set, and pull
back the resulting $D$-partitions to $\cK$.
The minimizers $\pi_tz$ remain defined over the original
compact set $\cK$, so all the preceding estimates still apply.
Proposition~\ref{prop:inputs}(\ref{item:contraction}) now gives
\[
 \gamma_2(\cK,D)
 \lesssim_d\sigma(A)+\frac1b+b[\sigma(A)+\beta(A)]^2.
\]
Choosing $b=[\sigma(A)+\beta(A)]^{-1}$ proves the claim.
\end{proof}

\subsection{Completion of the proof}

\begin{proof}[Proof of Theorem~\ref{thm:gaussian}]
The process $\ip{g}{\Phi(x)}$ is continuous on $\cK$,
has canonical pseudometric $D$, and vanishes at $x=0$.
Proposition~\ref{prop:inputs}(\ref{item:upperchaining}) and
Proposition~\ref{prop:transfer} give
\[
 \E\sup_{x\in\cK}|\ip{g}{\Phi(x)}|
 \lesssim\gamma_2(\cK,D)
 \lesssim_d\sigma(A)+\beta(A).
\]
For arbitrary $x^{(r)}\in B_2^{n_r}$, write
$x^{(r)}=x^{(r)}_+-x^{(r)}_-$, with both parts in
$B_{2,+}^{n_r}$. Expanding the multilinear form gives
\[
 \norm{(a_{\bi}g_{\bi})}_{\inj}
 \le2^d\sup_{x\in\cK}|\ip{g}{\Phi(x)}|,
\]
which proves the upper bound in \eqref{eq:coeff}.

For the lower bound, fixing all but one mode at coordinate
vectors gives
\[
 \norm{(a_{\bi}g_{\bi})}_{\inj}
 \ge\left(\sum_{\bi\in f}a_{\bi}^2g_{\bi}^2\right)^{1/2},
 \qquad f\in\cF.
\]
Jensen's inequality applied to the Euclidean norm of
$(a_{\bi}|g_{\bi}|)_{\bi\in f}$ yields
\[
 \E\left(\sum_{\bi\in f}a_{\bi}^2g_{\bi}^2\right)^{1/2}
 \ge\left(\sum_{\bi\in f}
              (a_{\bi}\E|g_{\bi}|)^2\right)^{1/2}
 =\sqrt{\frac2\pi}
   \left(\sum_{\bi\in f}a_{\bi}^2\right)^{1/2}.
\]
Maximizing over $f$ gives the lower bound in terms of $\sigma(A)$.
Testing individual entries and using Lemma~\ref{lem:gmax} also gives
\[
 \E\norm{(a_{\bi}g_{\bi})}_{\inj}
 \ge\E\max_{\bi\in\cI}a_{\bi}|g_{\bi}|
 \gtrsim\beta(A).
\]
Combining the two estimates proves the lower bound in
\eqref{eq:coeff}.
\end{proof}

\section{Sparse Gaussian tensors}\label{sec:masks}

\begingroup
In this section, we prove Theorem~\ref{thm:sparse} and
derive Corollaries~\ref{cor:simple-upper}
and~\ref{cor:bernoulli-tensors}.
\par\endgroup
For any nonnegative coefficient tensor
$B=(b_{\bi})_{\bi\in\cI}$, we use the fibre parameter $\sigma(B)$
and the entrywise parameter $\beta(B)$ introduced in
Sections~\ref{sec:main} and~\ref{sec:deterministic}, respectively.
The decreasing rearrangement in $\beta(B)$ includes zero entries.
Writing $A\odot\delta=(a_{\bi}\delta_{\bi})_{\bi\in\cI}$,
conditioning on the masks and applying
Theorem~\ref{thm:gaussian} reduces the proof of
Theorem~\ref{thm:sparse} to the estimate
\[
 \E_\delta\sigma(A\odot\delta)
 +\E_\delta\beta(A\odot\delta)
 \asymp_d \ext(A,P)+\vf(A,P).
\]

We first show that the entrywise term
$\E_\delta\beta(A\odot\delta)$ is comparable to $\ext(A,P)$
and establish the scalar truncation estimates.
We then average the fibre parameter in
Proposition~\ref{prop:averaging}.
\begingroup
The final subsection combines these estimates with the
conditional Gaussian bound and proves the two corollaries.
\par\endgroup

\subsection{The sparse entrywise maximum}\label{sec:extreme}

\begingroup
We now give an integral expression equivalent to the expected
maximum in \eqref{eq:expected-maximum} and relate it to the
ordered entrywise profile after masking.

\begin{lemma}[The sparse entrywise maximum]\label{lem:integral}
For every choice of profiles $A,P$,
\[
 \E_\delta\beta(A\odot\delta)\asymp\ext(A,P).
\]
Moreover, define
\begin{equation}\label{eq:extreme-profile}
 I(A,P)=\int_0^\infty
 \left[
   1\wedge
   \sum_{\substack{\bi\in\cI\\a_{\bi}>0}}
   p_{\bi}\exp\left(-\frac{t^2}{a_{\bi}^2}\right)
 \right]\,\dd t,
\end{equation}
where zero amplitudes are omitted from the sum. Then
\[
 \ext(A,P)\asymp I(A,P),
 \qquad
 \ext(A,P)\le\sqrt2\,I(A,P).
\]
\end{lemma}
\begin{proof}
Conditioning on the masks and applying Lemma~\ref{lem:gmax}
gives the first comparison.
Let $g$ be a standard Gaussian variable. Independence and
integration of the tail give
\[
 \ext(A,P)
 =\int_0^\infty
   \left[1-\prod_{\bi:\,a_{\bi}>0}
     \left(1-p_{\bi}\Prob\{|g|>t/a_{\bi}\}\right)\right]\dd t.
\]
For any finite family $z_i\in[0,1]$,
\[
 (1-\ee^{-1})\left(1\wedge\sum_i z_i\right)
 \le 1-\prod_i(1-z_i)
 \le 1\wedge\sum_i z_i.
\]
The lower bound follows from
$\prod_i(1-z_i)\le\exp(-\sum_i z_i)$, and the upper bound
is the union bound. The Gaussian tail satisfies
\[
 \exp(-u^2)\lesssim\Prob\{|g|>u\}
 \le\exp(-u^2/2),\qquad u\ge0.
\]
The lower estimate follows by integrating the density over
$[u,u+(u+1)^{-1}]$. For the upper estimate, shift the density
by $u$ and use
$\exp(-(s+u)^2/2)\le\exp(-u^2/2)\exp(-s^2/2)$ for $s\ge0$.
Substituting these bounds into the tail integral gives
$\ext(A,P)\gtrsim I(A,P)$ and
\[
 \ext(A,P)
 \le\int_0^\infty
 \left[1\wedge\sum_{\bi:\,a_{\bi}>0}
 p_{\bi}\exp\left(-\frac{t^2}{2a_{\bi}^2}\right)\right]\dd t
 =\sqrt2\,I(A,P).
\]
This proves the remaining assertions.
\end{proof}
\par\endgroup

We use $\ext(A,P)$ directly as the truncation scale.
The estimates below involving $A,P$ are immediate when
$\ext(A,P)=0$, since every entry with positive amplitude
then has activation probability zero. In their proofs,
we therefore assume $\ext(A,P)>0$.

\begin{lemma}[Contribution of large amplitudes]\label{lem:large}
For every choice of profiles $A,P$,
\[
 \sum_{\bi:\,a_{\bi}>4\ext(A,P)}p_{\bi}a_{\bi}
 \le4\ext(A,P).
\]
\end{lemma}
\begin{proof}
\begingroup
Conditional on the masks, test a position at which
$a_{\bi}\delta_{\bi}$ is largest. Using
$\E|g|=\sqrt{2/\pi}$ and \eqref{eq:expected-maximum} gives
\[
 \begin{aligned}
 \E_\delta\max_{\bi\in\cI}a_{\bi}\delta_{\bi}
 &\le\sqrt{\pi/2}\,
    \E_{\delta,g}\max_{\bi\in\cI}
       a_{\bi}\delta_{\bi}|g_{\bi}|\\
 &=\sqrt{\pi/2}\,\ext(A,P)
 \le2\ext(A,P).
 \end{aligned}
\]
\par\endgroup
By Markov's inequality, the probability that an entry with
amplitude greater than $4\ext(A,P)$ is activated is at most
$1/2$. Hence
\[
 \prod_{\bi:\,a_{\bi}>4\ext(A,P)}(1-p_{\bi})\ge\frac12.
\]
The events that exactly one of these entries is activated
are disjoint. Therefore
\[
 \begin{aligned}
 \E_\delta\max_{\bi\in\cI}a_{\bi}\delta_{\bi}
 &\ge
 \sum_{\bi:\,a_{\bi}>4\ext(A,P)}a_{\bi}p_{\bi}
 \prod_{\substack{\boldsymbol j:\,a_{\boldsymbol j}>4\ext(A,P)\\
                         \boldsymbol j\ne\bi}}
       (1-p_{\boldsymbol j})\\
 &\ge\frac12
 \sum_{\bi:\,a_{\bi}>4\ext(A,P)}a_{\bi}p_{\bi}.
 \end{aligned}
\]
Combining the two bounds proves the claim.
\end{proof}

The next estimate controls the second moment of the entrywise
profile after truncation. It will be used to bound the
fluctuations of the fibre sums.

\begin{lemma}[Second moment after truncation]\label{lem:beta-second}
Let $B=(b_{\bi})_{\bi\in\cI}$, where
$b_{\bi}=a_{\bi}\wedge4\ext(A,P)$. Then
\[
 \E_\delta\beta(B\odot\delta)^2\lesssim\ext(A,P)^2.
\]
\end{lemma}
\begin{proof}
Put $W=\max_{\bi\in\cI}b_{\bi}\delta_{\bi}|g_{\bi}|$.
\begingroup
Since $b_{\bi}\le a_{\bi}$, the definition
\eqref{eq:expected-maximum} gives $\E W\le\ext(A,P)$.
Markov's inequality and independence therefore imply
\par\endgroup
\[
 \prod_{\bi:\,b_{\bi}>0}
 \left(1-p_{\bi}\Prob\{|g|>4\ext(A,P)/b_{\bi}\}\right)
 \ge\frac12,
\]
where $g$ is a standard Gaussian variable. Using
$1-x\le\ee^{-x}$, we obtain
\[
 \sum_{\bi:\,b_{\bi}>0}
 p_{\bi}\Prob\{|g|>4\ext(A,P)/b_{\bi}\}\le\log2.
\]
For $v\ge u\ge0$, shifting the Gaussian density gives
\[
 \Prob\{|g|>v\}
 \le\exp\left(-\frac{v^2-u^2}{2}\right)
       \Prob\{|g|>u\}.
\]
Since $b_{\bi}\le4\ext(A,P)$, the preceding two estimates
and a union bound yield
\[
 \Prob\{W>t\}
 \le(\log2)\exp\left(
 -\frac{t^2-16\ext(A,P)^2}{32\ext(A,P)^2}\right),
 \qquad t\ge4\ext(A,P).
\]
Integration gives $\E W^2\lesssim\ext(A,P)^2$.
Conditional on the masks, Lemma~\ref{lem:gmax} and Jensen's
inequality give
\[
 \beta(B\odot\delta)^2
 \lesssim(\E_g W)^2\le\E_g W^2.
\]
Averaging over the masks proves the result.
\end{proof}

\subsection{Averaging the fibre profile}\label{sec:averaging}

To control the fluctuations of the fibre sums, we first
establish an estimate for their symmetrized counterparts.
The argument uses the fact that the fibres in each
$\cF_r$ are pairwise disjoint.

\begin{lemma}[Signed sums along fibres]\label{lem:signed-fibres}
For every deterministic nonnegative coefficient tensor $B$
and independent Rademacher variables $(\varepsilon_{\bi})$,
\[
 \E_\varepsilon\max_{f\in\cF}
 \left|\sum_{\bi\in f}\varepsilon_{\bi}b_{\bi}^2\right|
 \lesssim_d\beta(B)\sigma(B).
\]
\end{lemma}
\begin{proof}
The case $B=0$ is immediate. Fix a mode $r$. For each
$f\in\cF_r$, the Rademacher sum has sub-Gaussian variance
proxy
\[
 \sum_{\bi\in f}b_{\bi}^4
 \le\left(\max_{\bi\in f}b_{\bi}\right)^2\sigma(B)^2.
\]
Choose a largest coefficient in each fibre. Since the fibres
in $\cF_r$ are disjoint, these coefficients occupy distinct
positions. Their $j$th largest value is therefore at most
$b_j^*$. The sub-Gaussian tail bound and a union bound give,
for $t\ge2$,
\[
 \begin{aligned}
 &\Prob\left\{
 \max_{f\in\cF_r}
 \left|\sum_{\bi\in f}\varepsilon_{\bi}b_{\bi}^2\right|
 >t\beta(B)\sigma(B)\right\}\\
 &\qquad\le
 2\sum_{\substack{f\in\cF_r\\\max_{\bi\in f}b_{\bi}>0}}
 \exp\left(-\frac{t^2\beta(B)^2}
 {2(\max_{\bi\in f}b_{\bi})^2}\right)\\
 &\qquad\le2\sum_{j\ge1}(j+1)^{-t^2/2}
 \lesssim2^{-t^2/2}.
 \end{aligned}
\]
Integrating this bound gives the desired estimate for a
fixed mode, with an absolute comparison constant.
Since $\cF=\bigcup_{r=1}^d\cF_r$, summing over the modes
proves the stated bound.
\end{proof}

We can now average the fibre profile. The extreme term
$\ext(A,P)$ accounts for both the truncation error and the
fluctuations of the truncated fibre sums.

\begin{proposition}[Mask averaging]\label{prop:averaging}
For every choice of profiles $A,P$,
\[
 \ext(A,P)+\E_\delta\sigma(A\odot\delta)
 \asymp_d\ext(A,P)+\vf(A,P).
\]
\end{proposition}
\begin{proof}
Set $b_{\bi}=a_{\bi}\wedge4\ext(A,P)$ and
$B=(b_{\bi})_{\bi\in\cI}$. By the triangle inequality
and the bound of the maximum fibre norm by the entrywise
$\ell_1$ norm, Lemma~\ref{lem:large} gives
\[
 \E_\delta\sigma(A\odot\delta)
 \le\E_\delta\sigma(B\odot\delta)+4\ext(A,P).
\]
It remains to control the truncated tensor.
Recalling the definition of $v_s(A,P)$ in
\eqref{eq:truncated-variance}, symmetrization yields
\[
 \begin{aligned}
 \E_\delta\sigma(B\odot\delta)^2
 &\le v_{4\ext(A,P)}(A,P)^2
 +\E_\delta\max_{f\in\cF}
  \left|\sum_{\bi\in f}b_{\bi}^2(\delta_{\bi}-p_{\bi})\right|\\
 &\le v_{4\ext(A,P)}(A,P)^2
 +2\E_{\delta,\varepsilon}\max_{f\in\cF}
  \left|\sum_{\bi\in f}\varepsilon_{\bi}b_{\bi}^2\delta_{\bi}\right|.
 \end{aligned}
\]
For the second inequality, introduce an independent copy of
the masks, apply conditional Jensen's inequality, insert
independent signs in the differences, and use the triangle
inequality.
Conditional on the masks, Lemma~\ref{lem:signed-fibres}
applies to $B\odot\delta$. Together with
Lemma~\ref{lem:beta-second} and Cauchy--Schwarz, it gives
\[
 \begin{aligned}
 \E_\delta\sigma(B\odot\delta)^2
 &\lesssim_d v_{4\ext(A,P)}(A,P)^2
 +\E_\delta\bigl[
   \beta(B\odot\delta)\sigma(B\odot\delta)\bigr]\\
 &\lesssim_d v_{4\ext(A,P)}(A,P)^2
 +\ext(A,P)
   \bigl[\E_\delta\sigma(B\odot\delta)^2\bigr]^{1/2}.
 \end{aligned}
\]
Solving this quadratic inequality and using Jensen's
inequality, we obtain
\[
 \begin{aligned}
 \E_\delta\sigma(B\odot\delta)
 &\le\bigl[\E_\delta\sigma(B\odot\delta)^2\bigr]^{1/2}\\
 &\lesssim_d\ext(A,P)+v_{4\ext(A,P)}(A,P)
 \lesssim_d\ext(A,P)+\vf(A,P).
 \end{aligned}
\]
The last step uses
$a\wedge4\ext(A,P)\le4(a\wedge\ext(A,P))$.
Adding the truncation error proves the upper bound.

For the reverse bound, choose a fibre $f$ attaining the
maximum in the definition of $\vf(A,P)$, and put
\[
 Z=\sum_{\bi\in f}
       (a_{\bi}\wedge\ext(A,P))^2\delta_{\bi}.
\]
Independence gives
\[
 \E_\delta Z=\vf(A,P)^2,
\]
and
\[
 \operatorname{Var}(Z)
 =\sum_{\bi\in f}p_{\bi}(1-p_{\bi})
       (a_{\bi}\wedge\ext(A,P))^4
 \le\ext(A,P)^2\vf(A,P)^2.
\]
If $\vf(A,P)\ge\ext(A,P)$, then
$\E_\delta Z^2\le2\vf(A,P)^4$.
The Paley--Zygmund inequality therefore yields
\[
 \Prob\{Z\ge\vf(A,P)^2/2\}\ge\frac18,
 \qquad
 \E_\delta\sigma(A\odot\delta)
 \ge\E_\delta\sqrt Z\gtrsim\vf(A,P).
\]
If $\vf(A,P)<\ext(A,P)$, the extreme term already controls
$\vf(A,P)$. This proves the reverse bound.
\end{proof}

\subsection{Completion of the proofs}\label{sec:tensor-proofs}

The scalar and fibre estimates now allow us to average
the conditional Gaussian bound.

\begin{proof}[Proof of Theorem~\ref{thm:sparse}]
If $\ext(A,P)=0$, every masked entry vanishes almost surely,
and the conclusion is immediate. Otherwise, condition on
the masks and apply Theorem~\ref{thm:gaussian} to the
coefficient tensor $A\odot\delta$. Its comparison constants
depend only on $d$ and are uniform over all realizations
of the masks. Averaging gives
\[
 \begin{aligned}
 \E\norm{U}_{\inj}
 &\asymp_d
 \E_\delta\sigma(A\odot\delta)
 +\E_\delta\beta(A\odot\delta)\\
 &\asymp_d
 \E_\delta\sigma(A\odot\delta)+\ext(A,P)\\
 &\asymp_d\ext(A,P)+\vf(A,P).
 \end{aligned}
\]
The second comparison follows from Lemma~\ref{lem:integral},
and the last from Proposition~\ref{prop:averaging}.
All expectations are finite, since the index set $\cI$
is finite. This proves both bounds in \eqref{eq:main}.
\end{proof}

\begingroup
\begin{proof}[Proof of Corollary~\ref{cor:simple-upper}]
Let $N=\sum_{\bi:\,a_{\bi}>0}\delta_{\bi}$, so that $\E N=\Lambda$.
Conditional on the masks, Lemma~\ref{lem:gmax} gives
\[
 \E_g\max_{\bi\in\cI}a_{\bi}\delta_{\bi}|g_{\bi}|
 \lesssim a_*\sqrt{\log(1+N)},
\]
with both sides zero when $N=0$.
Since $x\mapsto\sqrt{\log(1+x)}$ is concave on $[0,\infty)$,
Jensen's inequality yields
\[
 \ext(A,P)\lesssim a_*\E\sqrt{\log(1+N)}
 \le a_*\sqrt{\log(1+\Lambda)}.
\]
Combining this with $\vf(A,P)\le V(A,P)$ and
Theorem~\ref{thm:sparse} proves the claim.
\end{proof}
\par\endgroup

We now deduce the Bernoulli estimate by symmetrization.

\begin{proof}[Proof of Corollary~\ref{cor:bernoulli-tensors}]
Let $\delta'$ be an independent copy of $\delta$, and recall
$q_{\bi}=2p_{\bi}(1-p_{\bi})$ and $Q=(q_{\bi})$.
The differences have the joint distribution
\[
 (\delta_{\bi}-\delta'_{\bi})_{\bi\in\cI}
 \overset{d}=(\varepsilon_{\bi}\eta_{\bi})_{\bi\in\cI},
\]
where the $\varepsilon_{\bi}$ are independent Rademacher variables,
the $\eta_{\bi}$ are independent Bernoulli variables with parameters
$q_{\bi}$, and the two families are independent.
Symmetrization and the Rademacher-to-Gaussian comparison give
\[
 \begin{aligned}
 \E\norm{A\odot(\delta-P)}_{\inj}
 &\le\E\norm{A\odot\varepsilon\odot\eta}_{\inj}\\
 &\le\sqrt{\pi/2}\,
       \E\norm{A\odot g\odot\eta}_{\inj}.
 \end{aligned}
\]
For the first inequality, use conditional Jensen's inequality
with respect to $\delta'$.
For the second, condition on $\eta$ and the signs, and apply
Jensen's inequality to the independent Gaussian magnitudes,
whose means are $\sqrt{2/\pi}$. The Gaussian signs are independent
of these magnitudes and have the same distribution as
$\varepsilon$.
Theorem~\ref{thm:sparse}, applied with probability profile $Q$,
now gives \eqref{eq:bernoulli-fluctuation}.
\end{proof}

\bigskip
\begingroup
\noindent\textbf{Acknowledgement.}\enspace
Ke Wang is supported by Hong Kong RGC grant GRF 16305526.
\par\endgroup

\bigskip
\begingroup
\noindent\textbf{Declaration of AI use.}\enspace
ChatGPT Astra was used to assist with mathematical discussions,
proof review, and language editing during the preparation of this paper.
The authors take full responsibility for the content and correctness
of the paper.
\par\endgroup

\end{document}